\documentclass[12pt]{article}
\usepackage{mathrsfs}
\usepackage{amsmath,amsfonts,amssymb,rotating,amsthm, bm}
\usepackage{hyperref}
\usepackage{array}
\usepackage{breqn}
\usepackage{color}
\usepackage{fancyheadings}
\usepackage{times}
\usepackage{enumitem}
\usepackage{dsfont}
\usepackage{caption}
\usepackage{subcaption}
\usepackage{booktabs}
\usepackage{tikz}
\usepackage{verbatim}
\usepackage{multirow}
\usepackage{appendix}
\usetikzlibrary{calc}
\usetikzlibrary{trees}

\tikzstyle{every node}=[circle,inner sep=1pt,fill=white!60]

\tikzstyle{tn}=[shape=circle, draw, color=black!70]
\tikzstyle{marke}=[shape=circle,minimum size=0.2cm, draw,blue]

\usepackage{times}    
\usepackage{mathptmx} 

\allowdisplaybreaks
\def\qed{\nopagebreak\hfill{\rule{4pt}{7pt}}}
\def\proof{\noindent {\it{Proof.} \hskip 2pt}}
\newtheorem{thm}{Theorem}[section]

\newtheorem{cor}[thm]{Corollary}

\theoremstyle{remark}

\numberwithin{equation}{section}

\usetikzlibrary{trees}

\tikzstyle{every node}=[circle,inner sep=1pt,fill=white!60]
\tikzstyle{tn}=[shape=circle, draw, color=black!70]
\tikzstyle{tn1}=[shape=circle, draw]
\tikzstyle{tn2}=[shape=circle, draw,inner sep=1.5pt]
\tikzstyle{tn5}=[shape=circle, inner sep=1.6pt,draw, color=black!70]
\tikzstyle{tn3}=[shape=rectangle, draw,inner sep=1.5pt]
\tikzstyle{tn4}=[shape=rectangle, draw,inner sep=1.5pt,color=black!70]
\tikzstyle{marke}=[shape=circle,minimum size=0.1cm, draw,blue]

\newcommand\gen{{\mathrm {Gen}}}

\newcommand\exc{{\mathrm{exc}}}
\newcommand\asc{\mathrm {asc}}
\newcommand\fix{\mathrm{fix}}

\newcommand\drop{\mathrm{drop}}
\newcommand\des{\mathrm{des}}
\newcommand\suc{\mathrm{suc}}
\newcommand\lsuc{\mathrm{lsuc}}
\newcommand\jump{\mathrm{jump}}

\begin{document}

\begin{center}
{\Large\bf
A Grammar of Dumont and

 a Theorem of  Diaconis-Evans-Graham
 }

\vskip 6mm

William Y.C. Chen$^1$ and Amy M. Fu$^2$

\vskip 3mm

$^{1}$Center for Applied Mathematics\\
Tianjin University\\
Tianjin 300072, P.R. China

\vskip 3mm

$^{2}$School of Mathematics\\
Shanghai University of Finance and Economics\\
Shanghai 200433, P.R. China

\vskip 3mm

Emails: { $^1$chenyc@tju.edu.cn, $^{2}$fu.mei@mail.shufe.edu.cn}

\vskip 6mm

\end{center}

\begin{abstract} We came across an unexpected connection between
a remarkable grammar of Dumont for the
joint distribution of $(\exc, \fix)$ over $S_n$ and
  a beautiful theorem of Diaconis-Evans-Graham
on  successions and fixed points of permutations.
With the grammar in hand, we
 demonstrate the advantage of the
grammatical calculus in deriving the generating functions,
where the constant property plays a substantial role.
On the grounds of
  left successions of a permutation,
 we present a grammatical treatment of  the joint
distribution investigated by Roselle.
Moreover, we obtain a left succession
    analogue
    of the Diaconis-Evans-Graham theorem, exemplifying the
    idea of a grammar assisted bijection. The grammatical
    labelings give rise to an equidistribution of
    $(\jump, \des)$ and $(\exc, \drop)$ restricted to the
    set of left successions
     and the set of fixed points,
    {where $\jump$ is  defined to be the number of ascents
    minus the number of left
    successions.}
        \end{abstract}



\section{Introduction}

This paper is concerned with a beautiful
theorem of Diaconis-Evans-Graham \cite{DEG-2014} on the
correspondence between successions and fixed
points of permutations. Unlike a typical equidistribution
property,
 an attractive feature  of this theorem
 is that the bijection can be restricted to permutations with a specific set of  {successions} and  permutations with the
same set of {fixed points}.

The topic of the enumeration of successions of permutations has a rich
history. Dumont referred to the work of Roselle \cite{Roselle-1968} on the
joint distribution of the number of ascents and the
number of successions. In fact, the
grammar proposed by Dumont \cite{Dumont-1996}
is meant to deal with the
joint distribution of the number of
excedances, the number of drops  and the number
of fixed points of a
permutation. His argument may be paraphrased in the
language of a grammatical labeling of complete increasing binary trees.
We will show that this grammar
is related to the Diaconis-Evans-Graham theorem, even though
it does not look so at
first sight.
It is worth mentioning that Dumont's citation to Roselle
was not accurate;
nevertheless, such an incident was somehow
just to the point. Indeed, this work would not
have come into being without the lucky pointer of Dumont.

First, we come to the realization that the grammar of Dumont can be adapted
 to a problem of Roselle. We just
 need to be more circumspect when
 it comes to
the notion of a left succession,
analogous to
that of a left peak of a permutation. In the
approach of Roselle, the consideration
of a left succession at position $1$
was considered informative for the computation
of the generating function of interior successions.
As for a left succession,
one assumes that a zero is patched at
 the beginning of a permutation.
  In contrast to a left succession, a usual
 succession is called an interior succession.

 Once the grammar is in place, a grammatical labeling
 is necessary in order to record a weighted counting of
 a combinatorial structure. A labeling scheme also
 makes it possible to carry out the grammatical calculus.
We will show how the grammar of Dumont
works for the joint distribution of $(\exc, \fix)$.
Furthermore, we give a different labeling scheme for
permutations which shows that the same grammar
of Dumont suits equally well  for the
joint distribution of $(\jump, \lsuc)$, where $\jump$ and $\lsuc$ denote the number of
jumps and the number of left successions of a permutation, respectively.
It is no surprise that
the constant property plays a substantial role in the
grammatical calculus.

While the grammar is instrumental in establishing an equidistribution,
it is not clear whether one can take a step forward in
obtaining a Diaconis-Evans-Graham type theorem
concerning a given set of left successions  and  the same set
fixed points. Fortunately, the answer is yes.
In fact, it is exactly
where the idea of a grammar assisted bijection comes on the scene.

\section{A grammar of Dumont }

In this section, we recall
a remarkable grammar of Dumont \cite{Dumont-1996}
for the joint distribution of the statistics
$(\exc, \drop, \fix)$ over $S_n$,   the set of
permutations of $[n]=\{1, 2,\ldots, n\}$, where $n\geq 1$.
For a permutation $\sigma = \sigma_1 \cdots \sigma_n \in S_n$,
an index $1 \leq i \leq n$ is called an excedance if
$\sigma_i > i$, or a drop if $\sigma_i < i$, or a fixed
point if $\sigma_i = i$. Clearly, $n$ cannot be an excedance and
$1$ cannot be a drop. The number of excedances, the
number of drops and the number of fixed points of
$\sigma$ are
denoted by $\exc(\sigma)$, $\drop(\sigma)$ and $\fix(\sigma)$,
respectively.  A drop of a permutation is also called an anti-excedance.

The joint distribution of $(\exc, \fix)$ was determined by Foata-Sch\"utzenberger
\cite{FS-1970}, see also Shin-Zeng \cite{SZ-2010}.
For $n\geq 1$,
define
\[ F_n(x,z) =\sum_{\sigma \in S_n}
  x^{\,\exc(\sigma)} z^{\,\fix(\sigma)} \]
and define $F_0(x,z)=1$.
 Then
\begin{equation} \label{fxz}
    \sum_{n=0}^\infty F_n(x,z) \frac{t^n}{n!}= \frac{(1-x)e^{zt}}{e^{xt}-xe^t}.
\end{equation}
Writing
\[ F_n(x,y,z) =\sum_{\sigma \in S_n}
  x^{\,\exc(\sigma)} y^{\,\drop(\sigma)} z^{\,\fix(\sigma)} \]
  and $F_0(x,y,z)=1$, \eqref{fxz}
can be converted into the homogeneous form
\begin{equation} \label{fsg-2}
\sum_{n=0}^\infty F_n(x,y,z) \frac{t^n}{n!}= \frac{(y-x)e^{zt}}{ye^{xt}-xe^{yt}}.
\end{equation}

Below are the first few values of $F_n(x,y,z)$:
\begin{eqnarray*}
    F_0(x,y,z) & = & 1,\\[3pt]
    F_1(x,y,z) & = & z,\\[3pt]
    F_2(x,y,z) & = & xy+z^2,\\[3pt]
    F_3(x,y,z) & = & 3xyz + x y^2 + x^2 y  +z^3,\\[3pt]
    F_4(x,y,z) & = & 6 x y z^2+4 x y^2 z+ x y^3 +4  x^2 y z+7 x^2 y^2+ x^3 y +z^4.
\end{eqnarray*}

The grammar of Dumont reads
\begin{equation} \label{GD}
   G=\{ a \rightarrow az, \; z\rightarrow xy, \;
     x\rightarrow xy, \; y \rightarrow xy\}.
\end{equation}

Let $D$ be the formal derivative with respect to $G$,
which can be expressed as a differential operator
\[ az {\partial \over \partial a}
+ xy {\partial \over \partial z}
+ xy {\partial \over \partial x}
+ xy {\partial \over \partial y}.\]
Dumont
\cite{Dumont-1996}
showed that
the polynomials $F_n(x,y,z)$ can be
generated by $D$.

\begin{thm}[Dumont]
    The following relation is valid for $n\geq 0$,
    \begin{equation}
    D^n(a) = a F_n(x,y,z).
    \end{equation}
\end{thm}

 Dumont's argument can be understood as a description of the
  procedure of recursively generating permutations
 in the cycle notation. Recall that a cycle is written
 in such a way that the minimum element is
 at the beginning and  the cycles of a permutation are arranged in the increasing order of the  minimum elements. Here we give an explanation in the
language of a grammatical labeling, which we call the
$(a,x,y,z)$-labeling, both for permutations and for
increasing binary trees.

Given a permutation $\sigma$ of $[n]$,   represent it
in the cycle notation. Use $a$ to signify the position where
a new cycle may be formed.
If $i$ is in a $1$-cycle, we label it by $z$.
If $(i,j)$ is an arc in the
cycle notation, that is, {$\sigma_i=j$},
we label it by $x$ if $i<j$, that is, $i$ is
an excedance, or by $y$ if $i>j$, that is, $i$ is a drop.
Then an insertion of $n+1$ into $\sigma$
can be formally described with the aid of the grammar rules.

For example,
below is  a permutation in the cycle notation, where the labels
are placed after each element  and the label $a$ is placed at the end:
\begin{equation} \label{ep}
(1\,x\,8\,y\,4\,x\,9\,y\,6 \, y) \ (2\,z )\ (3\,x\, 5\, y) \ (7 \, z)\; a.
\end{equation}

Relying on the grammar, one can build a complete increasing binary
tree to record the insertion process of generating a permutation of $[n+1]$
from a permutation of $[n]$, in the cycle notation, to be precise.
 To describe the procedure,
we represent a cycle by
arranging the minimum element at the
beginning followed by a permutation
of the remaining elements.
Clearly, this permutation following the minimum element
corresponds to a complete increasing binary tree, see,
for example, Stanley \cite[P. 23]{Stanley-2012}.

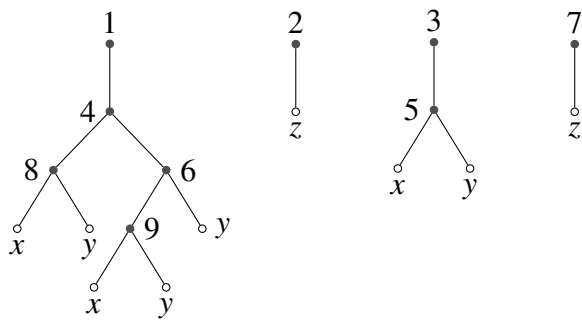
\begin{figure}[!ht]
\begin{center}

\begin{tikzpicture}[scale=0.6]
\node [tn,label=90:$1$]{}[grow=down]
    child {node [tn,label=180:{$4$}](four){}
       [sibling distance=25mm,level distance=13mm]
    child {node [tn,label=180:{$8$}](eight){}
     [sibling distance=16mm,level distance=13mm]
     child {node [tn1,label=-90:{$x$}](eightl){}}
     child {node [tn1,label=-90:{$y$}](eightr){}}
     }
      child {node [tn,label=0:{$6$}](six){}
     [sibling distance=16mm,level distance=13mm]
     child {node [tn,label=0:{$9$}](nine){}
     child {node [tn1, label=-90:{$x$}](ninel){}}
     child {node [tn1, label=-90:{$y$}](niner){}}
     }
     child {node [tn1,label=0:{$y$}](sixr){}}
     }
};
\end{tikzpicture}
\quad
\raisebox{68pt}{\begin{tikzpicture}[scale=0.6]
\node [tn,label=90:$2$]{}[grow=down]
    child {node [tn1,label=-90:{$z$}](four){}};
    \end{tikzpicture}}
\qquad
\raisebox{45pt}{\begin{tikzpicture}[scale=0.6]
\node [tn,label=90:$3$]{}[grow=down]
    child {node [tn,label=180:{$5$}](four){}
     [sibling distance=16mm,level distance=13mm]
    child {node [tn1,label=-90:{$x$}](eight){}}
    child {node [tn1,label=-90:{$y$}](eight){}}};
    \end{tikzpicture}}
    \qquad
\raisebox{68pt}{\begin{tikzpicture}[scale=0.6]
\node [tn,label=90:$7$]{}[grow=down]
    child {node [tn1,label=-90:{$z$}](four){}};
    \end{tikzpicture}}
\end{center}
\caption{A forest of planted increasing binary trees.}
\label{axyzf-2}
\end{figure}

Now, we may represent a cycle
by a planted complete increasing binary tree.
First, designate the minimum element as the root. If the cycle contains
only one element, then assign it a
$z$-leaf. Otherwise,  attach
the complete increasing binary tree corresponding to the permutation
following the minimum element as a subtree of the root.
Note that the external leaves of the complete increasing
binary tree comply with
the $(x,y)$-labeling for the
Eulerian polynomials, namely,
a left leaf is labeled by $x$ and a right
leaf is labeled by $y$. For example, the
cycles of the permutation in (\ref{ep}) are
represented by the forest of planted complete increasing
binary trees in
Figure \ref{axyzf-2}.

As the last step, we can put together these planted
increasing binary trees by drawing an edge between two roots next to each other to form a complete increasing binary tree
with the $(a,x,y,z)$-labeling for which
the root is $1$ and the rightmost leaf is labeled by $a$.

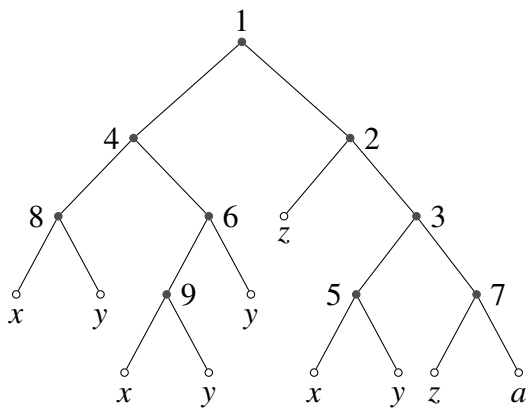
\begin{figure}[!ht]
\begin{center}
\begin{tikzpicture}[scale=0.8]
\node [tn,label=90:$1$]{}[grow=down]
	[sibling distance=36mm,level distance=16mm]
    child {node [tn,label=180:{$4$}](four){}
       [sibling distance=25mm,level distance=13mm]
    child {node [tn,label=180:{$8$}](eight){}
     [sibling distance=14mm,level distance=13mm]
     child {node [tn1,label=-90:{}](eightl){}}
     child {node [tn1,label=0:{}](eightr){}}
     }
      child {node [tn,label=0:{$6$}](six){}
     [sibling distance=14mm,level distance=13mm]
     child {node [tn,label=0:{$9$}](nine){}
     child {node [tn1, label=-90:{}](ninel){}}
     child {node [tn1, label=-90:{}](niner){}}
     }
     child {node [tn1,label=0:{}](sixr){}}
     }
     }
     child {node [tn,label=0:{$2$}](twoa){}
     [sibling distance=22mm,level distance=13mm]
      child {node [tn1,label=-90:{$z$}](twoal){}}
    child {node [tn,label=0:{$3$}](threea){}
       [sibling distance=20mm,level distance=13mm]
    child {node [tn,label=180:{$5$}](five){}
       [sibling distance=14mm,level distance=13mm]
      child {node [tn1,label=-90:{}](fivel){}}
    child {node [tn1,label=-90:{}](fiver){}}
    }
    child {node [tn,label=0:{$7$}](sevena){}
         [sibling distance=14mm,level distance=13mm]
    child {node [tn1,label=-90:{}](sevenl){}}
    child {node [tn1,label=-90:{$a$}](a){}}
    }
     }};
    \node [below=3pt] at (eightl){$x$};
     \node [below=3pt] at (eightr){$y$};
     \node [below=3pt] at (sixr){$y$};
     \node [below=3pt] at (ninel){$x$};
     \node [below=3pt] at (niner){$y$};
     \node [below=3pt] at (fivel){$x$};
     \node [below=3pt] at (fiver){$y$};
    \node [below=3pt] at (sevenl){$z$};
    \node [below=3pt] at (a){$a$};
\end{tikzpicture}
\end{center}
\caption{{The  $(a,x,y,z)$-labeling for $(\exc, \drop, \fix)$.}}
\label{f}
\end{figure}

For example, the forest in Figure \ref{axyzf-2}
can be put together into a complete increasing binary tree
 with an $(a,x,y,z)$-labeling
 in Figure \ref{f}.

We observe the
following properties.
\begin{itemize}
    \item
A $z$-leaf corresponds to a fixed point.
\item
An $x$-leaf corresponds to an excedance.

\item A $y$-leaf corresponds to a drop.
\end{itemize}

To recover a permutation $\sigma$
from a complete increasing binary tree $T$,
we may decompose $T$ into a forest of planted
increasing binary trees by removing the edges from the root
to the $a$-leaf and deleting the $a$-leaf.

The goal of this section is
to show that the generating function
of $F_n(x,y,z)$ can be easily derived by
the grammatical calculus.
A grammatical derivation of the
generating function of the Eulerian polynomials $A_n(x,y)$
was given in \cite{Chen-Fu-2022}. The same reasoning
can be carried over to the computation
of the generating function of $F_n(x,y,z)$.
Bear in mind that
the generating function with respect to the formal derivative $D$
is defined by
$$ \gen(w,t) = \sum_{n=0}^\infty D^n(w) \frac{t^n}{n!},$$
where  $w$ is a Laurent polynomial in the variables $a,x,y,z$.
Note that the generating function with respect to $D$ permits
the multiplicative property, which is equivalent to the Leibniz
rule, see \cite{Chen-Fu-2022} and references therein.

\begin{thm}\label{Thmat}
    We have
    \begin{equation} \label{gat}
        {\gen}(a, t) = \frac{a (y-x)e^{zt}}{y e^{xt}-xe^{yt}}.
    \end{equation}
\end{thm}

\noindent
{\it Proof.}
In virtue of  the rules
$$
x\rightarrow xy, \quad y \rightarrow xy,
$$
we obtain the generating function
$$
{\gen}(x,t)=\frac{x-y}{1-yx^{-1}e^{(x-y)t}},
$$
see \cite{Chen-Fu-2022}.
Since $D(z-y)=xy-xy=0$, i.e., $z-y$ is a constant relative to $D$,
we deduce that
$$
D^n(ax^{-1})=D^{n-1}\left(ax^{-1}(z-y)\right)=ax^{-1}(z-y)^n,
$$
and hence
\begin{equation}
    {\gen}(ax^{-1}, t) = \sum_{n=0}^\infty D^n(ax^{-1})\frac{t^n}{n!} = ax^{-1} e^{(z-y)t}.
\end{equation}

By the Leibniz rule or the product rule, we infer that
$$
{{\gen}(a,t)}={{\gen}(x\cdot ax^{-1},t )}={{\gen}(x,t)}\,{{\gen}(ax^{-1},t)}
=\frac{a (y-x)e^{zt}}{y e^{xt}-xe^{yt}},
$$
as required. \qed

Putting $a=1$, we arrive at Equation \eqref{fsg-2}.
Furthermore,  setting $z=0$ yields
the generating function of the derangement
polynomials; see Brenti \cite{Brenti-1990}.

\section{The  joint distribution of Roselle}

In this section, we give an account of
the generating function of Roselle \cite{Roselle-1968}
 for the joint distribution
of the number of ascents and the number of successions over $S_n$ in a nutshell.
Starting with recurrence relations, Roselle employed the
 symbolic method to accomplish the task of computation.
Such an antiquate
mechanism is rarely in demand these days, but perhaps it should
not be completely forgotten, even though it seems
obscure or dubious and even if its extinction might be inevitable.

\subsection{The formulas of Roselle}

Let us recall some definitions.
Let $n\geq 1$, and let $\sigma$ be a permutation of
$[n]$. We assume that $\sigma_0=0$. An ascent or a rise
of $\sigma$ is an index $0\leq i \leq n-1$ such that
$\sigma_i < \sigma_{i+1}$.
The number of ascents of $\sigma$ is denoted by $\asc(\sigma)$.
An index $i$ $(1\leq i \leq n-1)$ is called
a descent of $\sigma$ if $\sigma_i > \sigma_{i+1}$.
In this definition, the index $n$ is
not counted as a descent. The number of descents of
$\sigma$ is denoted by $\des(\sigma)$.
An index $i$ $(1\leq i \leq n-1)$ of $\sigma$ is
called a succession, or an interior succession,
if $\sigma_{i}+1 = \sigma_{i+1}$. We call an index $i$ $(1\leq i \leq n)$
a left succession if  $\sigma_{i-1}+1 = \sigma_{i}$.
Mind the subtlety  with respect to the range of indices
for a left succession.

In order to single out ascents that are not left successions,
we say that an index $1\leq i \leq n $ of $\sigma$ is
a jump if
$i-1$ is an ascent but $i$ is not a left succession, that is,
$\sigma_{i} \geq \sigma_{i-1} +2$.  The number of
jumps of $\sigma$ is denoted by $\jump(\sigma)$.

For $2\leq i \leq n$, if $i$ is a jump, then $i-1$
is called a big ascent by Ma-Qi-Yeh-Yeh \cite{MQYY-2024}, and the number of big ascents of $\sigma$ is denoted by
$\mathrm{basc}(\sigma)$. However, if $1$ is a jump,
it does not contribute to the counting of big ascents.

Let $P(n,r,s)$ denote the number
of permutations of
$[n]$ with $r$ ascents and $s$ (interior) successions.
For example, $P(3,2,0)=2$. The two permutations of
$\{1,2,3\}$ with
two ascents and no successions
are $132$, $213$. Nevertheless, $132$ has a left
succession. While the term of a left succession
is not manifestly put to use, one can find a clue through the
generating function for the
number of permutations of $[n]$ with $r$ ascents and
no left successions, see Roselle \cite{Roselle-1968}.

As to left successions, for $n\geq 1$,
let  $P^*(n,r,s)$ denote the number
of permutations of
$[n]$ with $r$ ascents and $s$ left successions.
Define $P_0^*(x,z)=1$ and define for $n\geq 1$,
$$
P_n^*(x,z)=\sum_{r=1}^{n}\sum_{s=0}^{r}P^*(n,r,s)x^{r-s}z^s.
$$
Since for $n\geq 1$ and for any permutation $\sigma\in S_n$,
\[ \asc(\sigma)=\jump(\sigma)+\lsuc(\sigma),
\]
we see that
$$
P_n^*(x,z)=\sum_{\sigma \in S_n}x^{\jump(\sigma)}z^{\lsuc(\sigma)}.
$$
The first few values of $P_n^*(x,z)$ are given below:
\begin{eqnarray*}
    P_0^* (x,z) & = & 1 ,\\[3pt]
    P_1^* (x,z) & = &  z,\\[3pt]
    P_2^* (x, z) & = &  x+z^2,\\[3pt]
    P_3^* (x, z) & = & x+x^2+3xz+z^3 ,\\[3pt]
    P_4^* (x, z) & = & x+7x^2+x^3+4xz+4x^2z+6xz^2+z^4.
\end{eqnarray*}

Below is the generating function of $P_n^*(x,z)$.

\begin{thm}[Roselle]
We have
\begin{equation}\label{gpstar}
 \sum_{n=0}^\infty P^{*}_n(x, z) \frac{t^n}{n!}
 =\frac{(1-x)e^{zt}}{e^{xt}-xe^t}.
\end{equation}
 \end{thm}

Notice that this formula coincides with \eqref{fxz}
for the joint distribution of
$(\exc,\fix)$. As will be seen,
this is by no means a coincidence.
We will encounter the same grammar in Section \ref{3.2}  and
so we ought to have the same
story.

Let us turn to the main theme of Roselle.
Set $P_0(x,z)=1$, and for $n\geq 1$ define
\begin{equation}
P_n(x,z)=\sum_{r=1}^n \sum_{s=0}^{r-1}P(n,r,s)x^rz^s,
\end{equation}
or equivalently,
\begin{equation}
  P_n(x,z)  =\sum_{\sigma \in S_n} x^{\,\asc(\sigma)}
  z^{\,\suc(\sigma)}.
\end{equation}
The first few values of $P_n (x,z)$ are given below:
\begin{eqnarray*}
    P_0  (x,z) & = & 1 ,\\[3pt]
    P_1  (x,z) & = &  x,\\[3pt]
    P_2  (x, z) & = & x^2z+x ,\\[3pt]
    P_3  (x, z) & = & x+2x^2+2x^2z+x^3z^2 ,\\[3pt]
    P_4  (x, z) & = & x+8x^2+2x^3+3x^2z+6x^3z+3x^3z^2+x^4z^3 .
\end{eqnarray*}
As shown by Roselle, the polynomials $P_n^*(x,z)$
serve as a stepping stone to compute $P_n(x,z)$.

\begin{thm}[Roselle] For $n\geq 1$, we have
\begin{equation} \label{pnstar}
P_n(x,z)=P^{*}_n(x,xz)+x(1-z)P^{*}_{n-1}(x,xz).
\end{equation}
\end{thm}

Combining the generating function of $P_n^*(x,z)$ and
the above relation gives rise to the
generating function of $P_n(x,z)$.

\begin{cor} We have
\begin{eqnarray} \label{gpn}
\sum_{n=0}^\infty P_{n+1}(x,z)\frac{t^n}{n!} =
\frac{x(1-x)^2e^{(xz+1)t}}{(e^{xt}-xe^t)^2}.
\end{eqnarray}
\end{cor}

\subsection{A grammatical labeling  for left successions} \label{3.2}

As alluded by the grammar of Dumont, we tend to believe that
the notion of a left succession should be considered
as a legitimate object of the subject, but
it does not seem to have gained enough recognition.

For $n \geq 1$, define
\[ L_n(x,y,z) = \sum_{\sigma \in S_n} x^{\,\jump(\sigma)}
y^{\,\des(\sigma)} z^{\,\lsuc(\sigma)}.\]
For $n=0$,  set $L_0(x,y,z)=1$.

 The following theorem shows that the polynomials $L_n(x,y,z)$
can be generated by the grammar $G$ of Dumont, that is,
$$
G=\{a \rightarrow az, \; z \rightarrow xy, \; x \rightarrow xy, \; y \rightarrow xy\} .
$$

\begin{thm}\label{ThmR}
Let $D$ be the formal derivative with respect to $G$.
For $n\geq 0$, we have
\begin{equation}
    D^n(a) = a L_n(x,y,z).
\end{equation}
\end{thm}

The above theorem can be justified by a labeling scheme of
permutations. Assume that $n\geq 1$ and
$\sigma$ is a permutation of $[n]$.
Consider the position after
each element $\sigma_i$ for $i=0,1,\ldots, n$, with $\sigma_0=0$.
First of all, label the position after the
maximum element $n$ by $a$.
{Next, if $\sigma_n \neq n$, label the position after $\sigma_n$ by $y$}. For the remaining
positions, if $i$ is a jump, then label the position on the {left} of $\sigma_i$ by $x$; if $i$ is a left succession, then label the position on the left of $\sigma_i$ by $z$, if $i$ is a descent and $\sigma_i \neq n$,
 label the position on the right of  $\sigma_i$ by $y$. Below is an example:
 \begin{equation}\label{ex-p9}
 0\ x\ 2 \ x \ 6 \ y \  3 \ z \ 4 \ y \
 1 \ x \  5 \  x \  8\  z\ 9 \ a \  7 \ y.
 \end{equation}

Write $*$ for the element $n+1$ to be
inserted into $\sigma$. The change of
labels can be described as follows. Assume that $*$ is to be inserted at the
position between $\sigma_{i}$ and $\sigma_{i+1}$, where
$0\leq i \leq n$. If $i=n$, $\sigma_{n+1}$ is considered empty.

\begin{enumerate}
    \item If $*$ is inserted at a position $a$, that is, $\sigma_{i}=n$,
   then we get $ n \,z\,  * \,a \sigma_{i+1}  $ in
    the neighborhood,
    this operation is captured by the rule $a \rightarrow az$.

    \item If $*$ is inserted at a position $x$,
   then we  see the update of $\sigma$:
   $\sigma_i \, x\, \sigma_{i+1} \rightarrow
   \sigma_i \, x \, * \, a\,  \sigma_{i+1}$. In the meantime,
    the label $a$ after $n$ in $\sigma$, wherever it is,
    will be switched to $y$, because $*$ is not inserted after $n$. This
    change of labels is reflected by the
    rule $x \rightarrow xy$.

    \item If $*$ is inserted at a position $y$, since $\sigma_i\not=n$,
    the update of $\sigma$ can be described by
    $\sigma_i \,
    y \, \sigma_{i+1} \rightarrow
    \sigma_{i}\, x \, * \, a \,\sigma_{i+1}$. In the
    meantime, the label $a$ after $n$
    in the labeling of $\sigma$, wherever it is,
    will be switched to $y$. This change of labels is governed by the rule
    $y\rightarrow xy$.

    \item If $*$ is inserted at a position $z$, then we have the update
    $\sigma_i \, z \sigma_{i+1} \rightarrow
    \sigma_i \, x \, * \, a \, \sigma_{i+1}$. In the meantime,
    the label  $a$ after $n$ in the labeling of $\sigma$, wherever it is,
    will be switched to $y$. This change of labels is in compliance
    with the rule $z\rightarrow xy$.

\end{enumerate}

We now have the same grammar for the two occasions. Thus
we are furnished with an equidistribution.

\begin{thm}
    For $n\geq 1$,
    the statistics $(\jump, \des, \lsuc)$  and the
    statistics $(\exc, \drop, \fix)$ are equidistributed
    over the set of permutations of $[n]$.
\end{thm}

In other words, the above theorem says that for $n\geq 0$,
\begin{equation} F_n(x,y,z)=L_n(x,y,z).
\end{equation}

In fact, we are going to pursue a stronger version
of the above theorem, that is, a left succession analogue of
the Diaconis-Evans-Graham theorem.
While a grammar might be sufficient to
guarantee an equidistribution
of two sets of statistics, it does not
tell us explicitly how to form a
bijection. Nevertheless, there are occasions  that
the grammar can be a guideline for establishing a
correspondence even under certain
constraints. We will come back to this point in Section 4.

\subsection{Back to interior successions}

Returning to the original formulation of the
joint distribution of Roselle, let $R_0(x,y,z)=1$,
and for $n\geq 1$, let
\begin{equation} \label{rnxyz}
R_n(x,y,z) = \sum_{\,\sigma\in S_n} x^{\,\jump(\sigma)} y^{\,\des(\sigma)} z^{\suc(\sigma)},
\end{equation}
which we call the Roselle polynomials. The first few values
of $R_n(x,y,z)$ are given below:
 \begin{eqnarray*}
    R_0(x,y,z) & = & 1,\\[3pt]
    R_1(x,y,z) & = & 1,\\[3pt]
    R_2(x,y,z) & = &  xy+z,\\[3pt]
    R_3(x,y,z) & = &  xy+2xyz+xy^2+ x^2y+z^2,\\[3pt]
    R_4(x,y,z) & = & 3 x y z
    +3 x y z^2 +x y^2+3 x y^2 z+ x y^3 +3 x^2 y z
    +x^2 y+7 x^2 y^2+ x^3 y+z^3 .
\end{eqnarray*}

Using the same reasoning for the grammatical labeling
for left successions together with a slight alternation of
the grammar, a grammatical calculus can be carried out for
the Roselle polynomials. Suppose that we are working with
the grammar for left successions, but
  we would like to avoid $1$ being counted as a left succession,
which is labeled by $z$. This requirement
can be easily met by turning to an additional label $b$ as a
substitute of the label $z$. That is to say, the rule
$z\rightarrow xy$  should be recast as
  $b\rightarrow xy$. For example,
we  should start with the initial
labeling $0b1a$ instead of $0z1a$.
As for the original labels $a,x,y,z$, their roles will
remain unchanged. Thus we meet with the mended
grammar:
\begin{equation} \label{gab}
    G=\{ a \rightarrow az,\; b\rightarrow  xy,   \;
     x \rightarrow xy, \; y\rightarrow xy, \; z\rightarrow xy\}.
\end{equation}
Let $D$ be the formal derivative of $G$ in \eqref{gab}.
Then we have
\[ D(ab) = abz + axy,\]
which is the sum of weights of the two permutations
\[0\  b \ 1 \ z \ 2 \ a , \quad 0\ x \  2 \ a \ 1 \ y. \]
In general, the polynomials $R_n(x,y,z)$ can also be
generated by the formal derivative $D$.

\begin{thm} \label{thm-R-n}
For $n\geq 1$, we have
\begin{equation}
    R_n(x,y,z) =  D^{n-1} (ab) \big| _{a=1,b=1}.
\end{equation}
\end{thm}

The grammatical calculus shows that the
generating function for the Roselle polynomials
is essentially a product of the generating function
of $L_n(x,y,z)$ and the generating function
of the bivariate Eulerian polynomials.

\begin{thm} We have
\begin{equation}\label{genab}
\gen(ab,t)=\frac{a (y-x)e^{zt}}{y e^{xt}-xe^{yt}} \left(\frac{x-y}{1-yx^{-1}e^{(x-y)t}}-x+b\right).
\end{equation}
   \end{thm}

\proof
By the Leibniz rule, we get
$$
\gen(ab,t)=\sum_{n=0}^\infty D^n(ab)\frac{t^n}{n!}=\gen(a,t)\gen(b,t).
$$
Since $D(b)=D(x)=xy$, it follows that
$$
\gen(b,t)=\gen(x,t)-x+b=\frac{x-y}{1-yx^{-1}e^{(x-y)t}}-x+b,
$$
which, together with Theorem \ref{Thmat}, implies
(\ref{genab}).  \qed

Next we show that the generating function
of $P_n(x,z)$ can be derived by using
the grammatical calculus.
Making substitutions in Theorem \ref{thm-R-n}  gives the following
relation.

\begin{cor} For $n\geq 1$, we have
    \begin{equation}
 P_n(x,z) = D^{n-1}(ab)|_{a=1,\, y=1, \,b=x,\, z=xz\,}.
\end{equation}
\end{cor}

\proof
Note that for any permutation $\sigma$ of $[n]$, we have for
$n\geq 1$,
\begin{equation}
    1+\jump(\sigma) + \suc(\sigma)= \asc(\sigma).
\end{equation}
By Theorem \ref{thm-R-n}, we find that
\begin{eqnarray*}
  D^{n-1}(ab)|_{a=1,\, y=1, \,b=x,\, z=xz\,} & = & x \sum_{\sigma\in S_n} x^{\,\jump(\sigma)} (xz)^{\,\suc(\sigma)} \\[6pt]
  & = &    \sum_{\sigma\in S_n} x^{\,\asc(\sigma)} z^{\,\suc(\sigma)},
\end{eqnarray*}
as required. \qed

The above relation enables us to
deduce the generating function of $P_n(x,z)$ from
that of $R_n(x,y,z)$, that is,
\begin{eqnarray*}
\sum_{n=0}^\infty P_{n+1}(x,z)\frac{t^n}{n!}
=\gen(ab,t)\big|_{a=1,\, y=1,\, b=x,\, z=xz}=\frac{x(1-x)^2e^{({x}z+1)t}}{(e^{xt}-xe^t)^2},
\end{eqnarray*}
which is in accordance with  \eqref{gpn}.

We finish this section with a relation between
$R_n(x,y,z)$ and $L_n(x,y,z)$, which can be
readily verified by the grammatical calculus.

\begin{thm} For $n\geq 0$, we have
\begin{equation} \label{Con-R-L}
    R_{n+1}(x,y,z) = L_{n}(x,y,z) + \sum_{k=1}^n
    \binom{n}{k} A_{k}(x,y) L_{n-k}(x,y,z),
\end{equation}
where for $k\geq 1$, $A_k(x,y)$  are the
bivariate Eulerian polynomials.
\end{thm}

This relation also admits a combinatorial interpretation.
Let $T$ be a complete increasing binary tree on $[n+1]$.
Suppose that we wish to interpret
$R_{n+1}(x,y,z)$ in terms of complete
increasing binary trees. We may adopt the following labeling
for $L_{n+1}(x,y,z)$, except that if the root of $T$ has a $z$-leaf,
we should label it by $1$ rather than $z$. If this is the case, then
the right subtree of $T$ can be viewed as a complete increasing tree on $[n]$
with a labeling, which contributes a term to
$L_n(x,y,z)$. If the root of $T$ has a
nonempty left subtree, then this left subtree does not have any $z$-leaves,
which can be reckoned as a labeling for the Eulerian polynomials, and
so we are through.

\section{An analogue
of the Diaconis-Evans-Graham theorem}

The main result of this paper is a left succession
analogue of the Diaconis-Evan-Graham theorem.
The grammar of Dumont can be
utilized to produce a bijection from permutations with a given set
of left successions to permutations with the same set
of fixed points, which possesses an additional equidistribution
property
concerning  $(\jump, \des)$
and $(\exc, \drop)$.

For $n\geq 1$ and a permutation $\sigma \in S_n$,  define
\begin{eqnarray*}
  M(\sigma) &  =  & \{ i \,
| \, 1 \leq i \leq n-1, \; \sigma_{i} + 1 = \sigma_{i+1} \},\\[3pt]
 G(\sigma) & = & \{ i \, | \,
 1 \leq i \leq n-1, \;
 \sigma_i = i \},\\[3pt]
 F(\sigma) & = & \{ i \, | \,
 1 \leq i \leq n , \;
 \sigma_i = i \}.
\end{eqnarray*}
It should be noted that the index $n$ is not
taken into consideration
in the definition of $G(\sigma)$.
Given a subset $I \subseteq [n-1]$, denote by
$M_n(I)$ the set of permutations of $[n]$ with $I$ being the set of (interior) successions, and denote by $G_n(I)$ the set of permutations $\sigma\in S_n$
such that $G(\sigma)=I$. Similarly, $F_n(I)$ denotes
the set of permutations $\sigma$ of $[n]$ such
that $F(\sigma)=I$.

\begin{thm}[Diaconis-Evans-Graham] \label{thm2.1}
Let $n\geq 1$ and $I\subseteq [n-1]$.
Then there is a bijection between
$M_n(I)$ and $G_n(I)$.
\end{thm}

For the special case $I=\emptyset$, a
permutation without successions is called a relative derangement.
Let $D_n$ denote the number of derangements of $[n]$, and
let $Q_n$ denote the number of relative derangements of $[n]$.
Roselle \cite{Roselle-1968} and
Brualdi \cite{Brualdi-2009} deduced that for $n \geq 1$,
\begin{equation} \label{qd}
    Q_n= D_n + D_{n-1}.
\end{equation}
A bijective proof of this relation was given in \cite{Chen-1996},
 appealing to the first fundamental transformation.
Taking $I=\emptyset$, a permutation in $G_n(I)$ may or may not have
$n$ as a fixed point. The permutations in these two cases are
counted by $D_{n-1}$ and $D_n$, respectively. Thus the $I=\emptyset$ case of
the proof of the Diaconis-Evans-Graham theorem
reduces to a combinatorial interpretation of (\ref{qd}).

Here comes the question of what happens for left successions.
To fit in the picture of
a grammar assisted bijection, we
find it more convenient to work with a variant or a reformulation
of the Diaconis-Evans-Graham theorem. Assume that
$n \geq 1$ and $\sigma\in S_n$.
Define
\begin{align*}
    \overline{M}(\sigma) &  =  \{ \sigma_i \, | \, 1 \leq i \leq n-1, \; \sigma_{i} + 1 = \sigma_{i+1} \}.
\end{align*}
It is readily seen that for any $\sigma\in S_n$,
\begin{eqnarray}
    \label{MMI}
    \overline{M}(\sigma^{-1}) &  =  & M(\sigma), \\
    G(\sigma^{-1}) & =  & G(\sigma),\\
    F(\sigma^{-1}) & = & F(\sigma).
\end{eqnarray}
where $\sigma^{-1}$ stands for the inverse of
$\sigma$.

Similar to the notation $\overline{M}(\sigma)$, for $n\geq 1$ and
a subset $I \subseteq [n-1]$, we define $\overline{M}_n(I)$
to be the set of permutations  $\sigma\in S_n$
such that $\overline{M}(\sigma) = I$.
Then
Theorem \ref{thm2.1} can be reformulated as follows.

\begin{thm}  \label{thm2.2}
Let $n\geq 1$ and $I\subseteq [n-1]$.
There is a bijection between
$\overline{M}_n(I)$ and $G_n(I)$.
\end{thm}

As a left succession analogue of $\overline{M}(\sigma)$, for $n\geq 1$ and a permutation $\sigma$ of $[n]$, we define
\[ \overline{L}(\sigma) = \{ \sigma_i \, | \, 1 \leq i \leq n, \; \sigma_{i-1} +1 = \sigma_i\}.\]
For a subset $I$ of $[n]$, define $\overline{L}_n(I)$ to
be the set of permutations $\sigma$ of $S_n$ such
that $\overline{L}(\sigma)=I$.

\begin{thm} \label{thm5.2}
For $n\geq 1$ and any $I \subseteq [n]$,
there is a bijection {  $\Phi$ from $\overline{L}_n(I)$ to ${F}_n(I)$}
that maps $(\jump, \des)$ to
  $(\exc, \drop)$.
\end{thm}
\proof
 Given a permutation $\sigma=\sigma_1 \cdots \sigma_n \in
 \overline{L}_n(I)$,
we wish to construct a complete increasing
binary $T$ with the $(a,x,y,z)$-labeling such that
  $\sigma_i \in \overline{L}(\sigma)$ if and only if
the vertex $\sigma_i$ has a $z$-leaf in $T$. Once the
correspondence is established, the equidistribution property
can be deduced from the interpretations of the labelings.

The map can be described as a recursive procedure.
For $n=1$, the permutation $z1a$ is mapped to the
complete increasing tree having one internal vertex
$1$  with a left $z$-leaf and a right $a$-leaf.

We now assume that $n\geq 1$ and  that $\sigma=\sigma_1\sigma_2\cdots \sigma_n $ is a permutation of $[n]$. As the induction
hypothesis, we assume that
$T$ is the tree corresponding to $\sigma$. For $1\leq i\leq n$,
the position $i$ is referred to   the position immediately before $\sigma_i$, whereas the position $n+1$ is meant to be the position after $\sigma_n$.

 To keep the procedure  running, we need to
  maintain additional properties of $\sigma$ and $T$.
 Besides having the same weight, they
 should be synchronized  in a certain sense.
 To be more specific,
 we say that the labeling of $\sigma$ is coherent with the labeling of $T$ provided that the following conditions are satisfied. In fact, these properties can be assured after each update.
 \begin{itemize}
  \item If the position $i$ in $\sigma$ is labeled by $x$, then the vertex $\sigma_{i}$ in $T$ has a $x$-leaf;
   \item If the position $i$ in $\sigma$ is labeled by $y$, then the vertex $\sigma_{i-1}+1$ in $T$ has a $y$-leaf;
      \item If the position $i$ in $\sigma$ is labeled by $z$, then the vertex $\sigma_{i}$ in $T$ has a $z$-leaf.
 \end{itemize}

Suppose that $*=n+1$ is to be inserted into $\sigma$.
It is necessary to find out how to update the tree $T$
accordingly.  Now that there are $n+1$ (insertion) positions of
$\sigma$ and there are $n+1$ leaves of $T$, it
suffices to  define a map from
the set of  positions
to the set of leaves of $T$ with the understanding that
when $*$ is inserted at a position, say $i$,
$T$ will be updated to $T'$ by
turning the corresponding leaf of $T$ into an internal vertex $*$. Denote by $\sigma'$ the permutation produced
from $\sigma$  by inserting $*$ at the position $i$.
There are four cases with regard to the four rules of the grammar.

\begin{enumerate}
    \item If $*$ is inserted at a position labeled by $a$, we add $*$ to $T$ at the position of the $a$-leaf. This operation is consistent with the rule $a \rightarrow az$.

   \item For a label $z$ at the position $i$,  by the induction hypothesis, we know that the vertex $\sigma_i$ in $T$
   has a $z$-leaf, so we can apply the rule $z\rightarrow xy$
    to this $z$-leaf to update $T$. Notice that when $*$ is inserted, the label $a$ on the right of $n$ in $\sigma$ will be switched to $y$. Observe that this $y$-label corresponds to the $y$-leaf of $*$ in $T'$. By inspection,
    we see that the labeling of $\sigma'$ is coherent with the labeling of $T'$.

    \item When the insertion occurs at position $i$ labeled by $x$,  by the induction hypothesis, we know that the vertex $\sigma_i$ in $T$ has a $x$-leaf. Then we apply the rule $x\rightarrow xy$ to this
    leaf. Notice that the $y$-leaf of $*$ in $T'$ corresponds to the $y$-label on the right of $n$.
    Again, it can be seen that
    the labeling of $\sigma'$ is coherent with the labeling of $T'$.

    \item  For a position $i$ labeled by $y$,   by the induction hypothesis, we know that the vertex $\sigma_{i-1}+1$ in $T$ has a $y$-leaf.  Then we can apply the rule $y \rightarrow xy$ to this leaf.    In this case, the labeling of $\sigma'$ remains coherent with the labeling of $T'$.

\end{enumerate}

So far we have provided a procedure to update
$T$ depending on where the element $*$ is inserted into $\sigma$.
Moreover, every stage of this procedure is reversible.
The detailed examination is omitted. As the grammar
ensures that the  map is weight-preserving,the weight of $\sigma$ equals that of $T$.

Observe that
a left succession, the element $\sigma_i$ for which
$\sigma_{i-1} +1   = \sigma_{i}$, to be
precise, is created in $\sigma$ whenever a vertex $\sigma_i$ with
a left $z$-leaf is created in $T$. Meanwhile,
a left succession $\sigma_i$ is destroyed in $\sigma$ whenever
a $z$-leaf with parent $\sigma_i$ is destroyed.

It should also be noted  that a jump value  $\sigma_i$ for which
$\sigma_{i-1}+2 \leq \sigma_{i}$ is created in $\sigma$ whenever a vertex $\sigma_i$ with
a left $x$-leaf is created in $T$. Meanwhile,
a jump value $\sigma_i$ is destroyed in $\sigma$ whenever
a left $x$-leaf with parent $\sigma_i$ is destroyed.

Since we have employed the cycle notation of a
permutation,  a vertex $\sigma_i$ with a left $z$-leaf
corresponds to a fixed point of a permutation, and
an $x$-leaf corresponds to an excedance, that is,
an element $\sigma_i$ such that $ \sigma_i > i$.
This completes the proof.
\qed

Figure \ref{ex} illustrates how to build
the corresponding trees step by step,
where an underlined label  indicates where an insertion takes place.

\begin{figure}[h!]
\begin{center}
\begin{tabular}{lc}
Permutations & Trees\\[5pt]
\raisebox{15pt}{$z\ 1\ {\underline a}$} &  \begin{tikzpicture}[scale=0.6]
\node [tn,label=90:{$1$}]{}[grow=down]
	[sibling distance=16mm,level distance=10mm]
     child {node [tn1,label=-90:{$z$}](){}}
     child {node [tn1,label=-90:{${\underline a}$}](){}};
\end{tikzpicture}\\[5pt]
\raisebox{30pt}{$z\ 1\ {\underline   z}
\ 2\ a$} & \begin{tikzpicture}[scale=0.6]
\node [tn,label=90:{$1$}]{}[grow=down]
	[sibling distance=16mm,level distance=10mm]
     child {node [tn1,label=-90:{$z$}](){}}
     child {node [tn,label=0:{$2$}](){}
     child {node [tn1,label=-90:{${\underline z}$}](){}}
     child {node [tn1,label=-90:{$a$}](){}} };
\end{tikzpicture}\\[5pt]
\raisebox{40pt}{$z\ 1\ x\ 3\ a\ 2\ {\underline y}$} & \begin{tikzpicture}[scale=0.6]
\node [tn,label=90:{$1$}]{}[grow=down]
	[sibling distance=16mm,level distance=10mm]
     child {node [tn1,label=-90:{$z$}](){}}
     child {node [tn,label=0:{$2$}](){}
     child {node [tn,label=180:{$3$}](){}
     child {node [tn1,label=-90:{$x$}](){}}
     child {node [tn1,label=-90:{${\underline y}$}](){}}}
     child {node [tn1,label=-90:{$a$}](){}} };
\end{tikzpicture}\\[5pt]
\raisebox{45pt}{$z\ 1\ x\ 3\ y\ 2\ x\ 4\ {\underline a}$} & \begin{tikzpicture}[scale=0.6]
\node [tn,label=90:{$1$}]{}[grow=down]
	[sibling distance=16mm,level distance=10mm]
     child {node [tn1,label=-90:{$z$}](){}}
     child {node [tn,label=0:{$2$}](){}
     child {node [tn,label=180:{$3$}](){}
     child {node [tn1,label=-90:{$x$}](){}}
     child {node [tn,label=180:{$4$}](){}
    child {node [tn1,label=-90:{$x$}](){}}
     child {node [tn1,label=-90:{$y$}](){}}}}
     child {node [tn1,label=-90:{$\underline{a}$}](){}} };
\end{tikzpicture}\\[5pt]
\raisebox{45pt}{$z\ 1\ \underline{x}\ 3\ y\ 2\ x\ 4\ z\ 5\ a$} &
\begin{tikzpicture}[scale=0.6]
\node [tn,label=90:{$1$}]{}[grow=down]
    [sibling distance=22mm,level distance=10mm]
    child {node [tn1,label=-90:{$z$}]{}}
    child {node [tn,label=0:{$2$}]{}
    [sibling distance=22mm,level distance=10mm]
    {
        child {node [tn,label=180:{$3$}]{}
        [sibling distance=15mm,level distance=10mm]{
            child {node [tn1,label=-90:{$\underline{x}$}]{}}
            child {node [tn,label=180:{$4$}]{}
            [sibling distance=11mm,level distance=10mm]{
                child {node [tn1,label=-90:{$x$}]{}}
                child {node [tn1,label=-90:{$y$}]{}}
            }}
        }}
        child {node [tn,label=0:{$5$}]{}
        [sibling distance=11mm,level distance=10mm]{
            child {node [tn1,label=-90:{$z$}]{}}
            child {node [tn1,label=-90:{$a$}]{}}
        }}
    }};
\end{tikzpicture}\\[5pt]

\raisebox{45pt}{$z\ 1\ x\ 6\ a\ 3\ y\ 2\ x\ 4\ z\ 5\ y$} \qquad \qquad &
\begin{tikzpicture}[scale=0.6]
\node [tn,label=90:{$1$}]{}[grow=down]
    [sibling distance=23mm,level distance=10mm]
    child {node [tn1,label=-90:{$z$}]{}}
    child {node [tn,label=0:{$2$}]{}[sibling distance=23mm,level distance=10mm]{
        child {node [tn,label=180:{$3$}]{}
        [sibling distance=18mm,level distance=10mm]{
            child {node [tn,label=180:{$6$}]{}
            [sibling distance=11mm,level distance=10mm]{
                child {node [tn1,label=-90:{$x$}]{}}
                child {node [tn1,label=-90:{$y$}]{}}
            }}
            child {node [tn,label=180:{$4$}]{}
            [sibling distance=11mm,level distance=10mm]{
                child {node [tn1,label=-90:{$x$}]{}}
                child {node [tn1,label=-90:{$y$}]{}}
            }}
        }}
        child {node [tn,label=0:{$5$}]{}
        [sibling distance=11mm,level distance=10mm]{
            child {node [tn1,label=-90:{$z$}]{}}
            child {node [tn1,label=-90:{$a$}]{}}
        }}
    }};
\end{tikzpicture}\\[5pt]
\end{tabular}
\end{center}
\caption{An example.}
\label{ex}
\end{figure}

For $n=3$, the correspondence is given in
the table below.  The
cases when $\overline{L}_n(I)=\emptyset$ or $F_n(I)=\emptyset$  are not listed, such as $I=\{1,2\}$.

\begin{center}
\renewcommand\arraystretch{1.33}
\begin{tabular}{|c|c|c|c|}
\hline
$I\subseteq [n]$ & $\overline{L}_n(I)$ & ${F}_n(I)$ & $(\text{jump}, \text{des})~\text{of}~\overline{L}_n(I) \leftrightarrow (\text{exc},\text{drop})~\text{of}~F_n(I)$ \\
\hline
\multirow{2}{*}{$\emptyset$} & $2\;1\;3$ & $(1\,2\,3)$ & $(2,1)$ \\
\cline{2-4}
& $3\;2\;1$  & $(1\,3\,2)$ & $(1,2)$ \\
\hline
$\{1\}$ & $1\;3\;2$  & $(1)(2\,3)$ & $(1,1)$ \\
\hline
$\{2\}$ & $3\;1\;2$ & $(1\,3)(2)$ & $(1,1)$ \\
\hline
$\{3\}$  & $2\;3\;1$ & $(1\,2)(3)$ & $(1,1)$ \\
\hline
$\{1,2,3\}$ & $1\;2\;3$ & $(1)(2)(3)$ & $(0,0)$ \\
\hline
\end{tabular}
\end{center}

To conclude, we remark that the above grammar
assisted bijection permits a refined equidistribution
property in terms of set-valued statistics. As shown in
\cite{MQYY-2024}, a grammar may be a helpful platform
to deal with set-valued statistics.
Roughly speaking, the above grammar assisted bijection maps
elements associated with the $x$-labels in a permutation
to elements associated with the $x$-labels in a complete
increasing binary tree.
More precisely, let
\begin{eqnarray*}
    \overline{\rm Jump}(\sigma) & = & \{\sigma_i \, | \, 1 \leq i \leq n, \;
        \sigma_{i-1} + 2 \leq \sigma_i \}, \\[3pt]
    \overline{\rm Exc}(\sigma) & = & \{\sigma_i \, | \, 1 \leq i \leq n, \;
        \sigma_{i } > i \}.
    \end{eqnarray*}
In other words, the set $\overline{\rm Jump}(\sigma)$ consists
of elements immediately to the right of the $x$-labels of $\sigma$, whereas
the elements in $\overline{\rm Exc}(\sigma)$ are exactly the
vertices having an $x$-leaf in a complete increasing binary tree.
Thus for the bijection $\Phi$ in the theorem and for any permutation
$\sigma$ of $S_n$, we have
\begin{equation}  \overline{\rm Jump}(\sigma)  =
 \overline{\rm Exc}(\Phi(\sigma)).
 \end{equation}
 For example, let $\sigma = 1\,6\,3\,2\,4\,5$.
 Then we have $\Phi(\sigma)=(1)(2\,6\,3\,4)(5)$. It is
 readily checked that
\[         \overline{L}(\sigma) =  {F}(\Phi(\sigma))=\{1,5\} \]
and \[ \overline{\rm Jump}(\sigma)
=\overline{\rm Exc}(\Phi(\sigma))=\{4,6\}.
\]
 Similarly, the $y$-labels are
 related to the set-valued refinements of
 $\des$ and $\drop$. So our grammar assisted
 bijection suits the purpose of producing a
   set-valued equidistribution.

\vskip 6mm \noindent{\large\bf Acknowledgments.} We are grateful to the
referee for insightful comments and
substantial suggestions.
 This work was supported
by the National
Science Foundation of China.


\begin{thebibliography}{99}



\bibitem{Brenti-1990}
F. Brenti, Unimodal polynomials arising from symmetric functions,
Proc. Amer. Math. Soc., 108  (1990) 1133--1141.


\bibitem{Brualdi-2009}
R.A. Brualdi, Introductory Combinatorics, 5th ed., Pearson/Prentice Hall, 2009.


 \bibitem{Chen-1996}
W.Y.C. Chen, The skew, relative, and classical derangements,
Discrete Math., 160 (1996) 235--239.


\bibitem{Chen-Fu-2022}
W.Y.C. Chen and A.M. Fu,
A Context-free grammar for the $e$-positivity
of the trivariate second-order Eulerian polynomials,
Discrete Math., 345  (2022) 112661.



\bibitem{DEG-2014}
P. Diaconis, S.N. Evans and R. Graham,
Unseparated pairs and fixed points in random permutations,
Adv. in Appl. Math., 61 (2014) 102--124.


\bibitem{Dumont-1996}
D. Dumont, Grammaires de William Chen et d\'erivations dans les arbres et arborescences, S\'em. Lothar. Combin., 37 (1996) Art. B37a.


 \bibitem{FS-1970}
 D. Foata and M.-P. Sch\"utzenberger,
Th\'eorie g\'eom\'etrique des polyn\^omes eul\'eriens,
Lecture Notes in Math., Vol. 138, Springer-Verlag, Berlin-New York, 1970.


\bibitem{MQYY-2024}
S.-M. Ma, H. Qi, J. Yeh and Y.-N. Yeh,
On the joint distributions of succession and Eulerian statistics,
arXiv:2401.01760.


\bibitem{Roselle-1968}
D.P.  Roselle,
Permutations by number of rises and successions,
Proc. Amer. Math. Soc.,
19  (1968) 8--16.

\bibitem{SZ-2010}
H. Shin and J. Zeng,
The $q$-tangent and $q$-secant numbers via continued fractions,
 European J. Combin., 31  (2010) 1689--1705.

\bibitem{Stanley-2012}
R.P. Stanley, Enumerative Combinatorics, Vol. I, second ed.,
Cambridge Univ. Press, Cambridge, 2012.

\end{thebibliography}
\end{document}